\documentclass [11pt]{article}
\usepackage{mathrsfs}
\usepackage{amsthm,amsmath,amsfonts,color}
\usepackage{latexsym}

\newtheorem{thm}{Theorem}

\newtheorem{cor}{Corollary}

\newtheorem{prob}{Problem}
\newtheorem{prop}{Proposition}

\def\qed{\hfill \rule{2.5mm}{2.5mm}}
\date{}

\begin{document}

\title{Bounds on the 2-rainbow domination number of graphs}

\author{Yunjian Wu$^{1}$\thanks{Corresponding author: y.wu@seu.edu.cn}\ \ and N. Jafari Rad$^{2}$
\\ {\small $^1$ Department of Mathematics}
\\ {\small \small Southeast University, Nanjing 211189, China}
\\ {\small $^2$ Department of Mathematics}
\\ {\small Shahrood University of Technology, Shahrood, Iran}
}

\maketitle

\begin{abstract}
A {\it $2$-rainbow domination function} of a graph $G$ is a function
$f$ that assigns to each vertex a set of colors chosen from the set
$\{1,2\}$, such that for any $v\in V(G)$, $f(v)=\emptyset$ implies
$\bigcup_{u\in N(v)}f(u)=\{1,2\}$. The {\it $2$-rainbow domination
number $\gamma_{r2}(G)$} of a graph $G$ is the minimum
$w(f)=\Sigma_{v\in V}|f(v)|$ over all such functions $f$. Let $G$ be
a connected graph of order $|V(G)|=n\geq 3$. We prove that
$\gamma_{r2}(G)\leq 3n/4$ and we characterize the graphs achieving
equality. We also prove a lower bound for 2-rainbow domination
number of a tree using its domination number. Some other lower and
upper bounds of $\gamma_{r2}(G)$ in terms of diameter are also
given.

\bigskip
\begin{flushleft}
\end{flushleft}


\end{abstract}

\vskip 1mm \vskip 1mm

\section{Introduction}

We follow the notation of \cite{bb} in this paper. Specifically, let $G = (V, E)$
be a graph with vertex set $V$ and edge set $E$. $P_k$ and $C_k$ denote
a path and a cycle of order $k$, respectively. For any vertex $v\in
V$, the {\it open neighborhood of $v$} is the set $N(v)=\{u\in V\ |\
uv\in E\}$ and the {\it closed neighborhood} is the set
$N[v]=N(v)\cup\{v\}$. For a set $S\subseteq V$, the open
neighborhood is $N(S)=\bigcup_{v\in S}N(v)$ and the closed
neighborhood is $N[S]=N(S)\cup S$. The {\it diameter} of $G$ is the
maximum distance between vertices of $G$, denoted by $diam(G)$. A
{\it penultimate vertex} is any neighbor of a vertex with degree one
(the vertex of degree one is also called a {\it leaf} in a tree), and
a {\it pendent} edge is an edge incident with a vertex of degree
one. A {\it star} is a tree isomorphic to a bipartite graph
$K_{1,k}$ for $k\geq 1$. A {\it double-star $DS_{r,s}$} is a tree
with diameter 3 in which there are exactly two penultimate vertices
with degrees $r+1$ and $s+1$, respectively. A set $S\subseteq V$ is
a {\it dominating set} of $G$ if every vertex not in $S$ is adjacent
to a vertex in $S$. The {\it domination number} of $G$, denoted by
$\gamma(G)$, is the minimum cardinality of a dominating set. A
thorough study of domination concepts appears in \cite{hay}. For a
pair of graphs $G$ and $H$, the {\it Cartesian product} $G\square H$ of
$G$ and $H$ is the graph with vertex set $V(G)\times V(H)$, where
two vertices are adjacent if and only if they are equal in one
coordinate and adjacent in the other.

Let $f$ be a function that assigns to each vertex a set of colors
chosen from the set $\{1,\cdots, k\}$; that is, $f: V(G)\rightarrow
\mathscr{P}(\{1,\cdots, k\})$. If for each vertex $v\in V(G)$ such
that $f(v)=\emptyset$. we have
$$\bigcup_{u\in N(v)}f(u)=\{1,\cdots, k\}.$$
Then $f$ is called a {\it $k$-rainbow dominating function} ($k$RDF)
of $G$. The weight, $w(f)$, of a function $f$ is defined as
$w(f)=\Sigma_{v\in V(G)}|f(v)|$. The minimum weight of a $k$-rainbow
dominating function is called the {\it $k$-rainbow domination
number} of $G$, which we denote by $\gamma_{rk}(G)$. We say that a
function $f$ is a {\it $\gamma_{rk}(G)$-function} if it is a $k$RDF
and $w(f)=\gamma_{rk}(G)$. The concept of rainbow domination was
introduced in \cite{bhr}, and used in obtaining some bounds on the
paired-domination number of Cartesian products of graphs, see also
\cite{bhr0}. A more ambitious motivation for the introduction of
this invariant was inspired by the following famous open problem
\cite{Vizing}:

\medskip

\noindent{\bf Vizing's Conjecture}. For any graphs $G$ and $H$,
$\gamma(G\square H) \geq \gamma(G)\gamma(H)$.

\medskip

In the language of domination of Cartesian products,
Hartnell and Rall \cite{har} obtained a couple of observations about
rainbow domination, for instance, $\min\{|V(G)|,\gamma(G)+k-2\}\leq
\gamma_{rk}(G)\leq k\gamma(G)$. Rainbow domination of a graph $G$ coincides with the ordinary
domination of the Cartesian product of $G$ with the complete graph,
in particular $\gamma_{r2}(G)=\gamma(G\square K_2)$ for any graph
$G$ \cite{bhr}. Notably a lower bound for the
2-rainbow domination number of a graph expressed in terms of its
ordinary domination could bring a new approach to the much desired
proof of Vizing's conjecture. In particular, Bre\v{s}ar, Henning and Rall \cite{bhr} proposed the following problem:

\begin{prob}\label{b}
{\rm{(Bre\v{s}ar, Henning and Rall \cite{bhr})}}. For any graphs $G$ and $H$,
$\gamma_{r2}(G\square H) \geq \gamma(G)\gamma(H)$.
\end{prob}

A {\it Roman domination function} of a graph $G$ is a function $g:
V\rightarrow\{0,1,2\}$ such that every vertex with 0 has a neighbor
with 2. The {\it Roman domination number} $\gamma_R(G)$ is the
minimum of $g(V(G))=\Sigma_{v\in V}g(v)$ over all such functions. In \cite{wx}, the authors showed the following result:

\begin{thm}\label{w}
{\rm{(Wu and Xing \cite{wx})}}Let $G$ be a graph. Then $\gamma(G)\leq \gamma_{r2}(G)\leq \gamma_R(G)\leq 2\gamma(G).$
\end{thm}

In \cite{wu}, Wu showed the following weaker form of Problem \ref{b} by Theorem \ref{w}:

\begin{thm}\label{wu}
{\rm{(Wu \cite{wu})}}For any graphs $G$ and $H$,
$\gamma_R(G\square H) \geq \gamma(G)\gamma(H).$
\end{thm}

In fact, Both Problem \ref{b} and Theorem \ref{wu} are improvements of the result given by Clark and Suen \cite{CS}:

\begin{thm}\label{c}
{\rm{(Clark and Suen \cite{CS})}} For any graphs $G$ and $H$, $2\gamma(G\square H)
\geq \gamma(G)\gamma(H)$.
\end{thm}

Nevertheless the concept of rainbow domination seems to be of
independent interest as well and it attracted several authors who
provided structural and algorithmic results on this invariant
\cite{bs,cwz, tlyl,xu}. In particular, it was shown that the problem
of deciding if a graph has a 2-rainbow dominating function of a
given weight is NP-complete even when restricted to bipartite graphs
or chordal graphs \cite{bs}. Also a few exact values and bounds for
the 2-rainbow domination number were given for some special classes
of graphs, including generalized Petersen graphs \cite{bs,xu}.

For a graph $G$, let $f: V(G)\rightarrow \mathscr{P}(\{1,2\})$ be a
2RDF of $G$ and  $(V_0, V_1^1, V_1^2, V_2)$ be the ordered partition
of $V(G)$ induced by $f$, where $V_0=\{v\in V(G)\ |\
f(v)=\emptyset\}$, $V_1^1=\{v\in V(G)\ |\ f(v)=\{1\}\}$,
$V_1^2=\{v\in V(G)\ |\ f(v)=\{2\}\}$ and $V_2=\{v\in V(G)\ |\
f(v)=\{1,2\}\}$. Note that there exists a 1-1 correspondence between
the functions $f: V(G)\rightarrow \mathscr{P}(\{1,2\})$ and the
ordered partitions $(V_0, V_1^1, V_1^2, V_2)$ of $V(G)$. Thus we
will write $f=(V_0, V_1^1, V_1^2, V_2)$ for simplicity.

In this paper we present some general bounds on the 2-rainbow
domination number of a graph that are expressed in terms of the
order and domination number of a graph. More specifically, we show
that $\gamma_{r2}(G)\leq 3|V(G)|/4$ and we characterize the graphs
achieving equality. We also prove a lower bound for the 2-rainbow domination number of a tree using its domination number. The latter
lower bound goes in the direction of the original goal, mentioned
above, to obtain a new approach for establishing Vizing's
conjecture. Some other lower and upper bounds of $\gamma_{r2}(G)$ in
terms of diameter are also given.

\medskip
\section{Main results}

Our aim in this section is to determine some bounds on the 2-rainbow domination number of graphs.

\subsection{Upper bounds}

We first recall a few definitions. A {\it subdivision} of an edge
$uv$ is obtained by removing edge $uv$, adding a new vertex $w$, and
adding edges $uw$ and $vw$. Let $t\geq 2$. A {\it spider} $(${\it
wounded spider}$)$ is the graph formed by subdividing some edges (at
most $t-1$ edges) of a star $K_{1, t}$. The unique center of $K_{1,t}$ is
also called the {\it center} of the spider. Only one vertex of the spider $P_4$ can be
called the center.

\begin{prop}\label{spider}
Let $G$ be a spider of order $|V(G)|=n\geq 3$, then
$\gamma_{r2}(G)\leq 3n/4$. Moreover, the equality only holds for a
path of order four.
\end{prop}

\proof Let $u$ be the center of $G$. Suppose $u$ has $x$ penultimate
neighbors and $y$ non-penultimate neighbors. Then $n=2x+y+1$.

If $x\geq 3$ or $y\geq 2$, we set

$$f(v)=\left\{\begin{array}{lllcrrr}
    \{1,2\}  & \mbox{\ \ $v=u$},\\[5pt]
    \{1\}\ \mbox{or}\ \{2\} & \mbox{\ \ $v \mbox{\ is\ at\ distance\ two to $u$}$},\\[5pt]
    \hskip0,2cm \emptyset  & \mbox{\ \ otherwise}.\\[5pt]

\end{array} \right.$$

If $x=2$ and $y\leq 1$, we set

$$f(v)=\left\{\begin{array}{lllcrrr}
    \{1\}  & \mbox{\ \ $v=u$},\\[5pt]
    \{2\}  &  \ \ \mbox{$v$\ is\ a\ leaf},\\[5pt]
     \hskip0,2cm \emptyset  & \mbox{\ \ otherwise}.\\[5pt]

\end{array} \right.$$

In both cases, $\gamma_{r2}(G)\leq w(f)< 3n/4$.

If $x=y=1$, then $G$ is a path of order four. Clearly,
$\gamma_{r2}(G)=3=3n/4$. \qed

\medskip

\begin{thm}\label{tree}
Let $T$ be a tree of order $n\geq 3$, then
$\gamma_{r2}(T)\leq 3n/4$.
\end{thm}

\proof We use induction on $n$. The base step handles trees with few
vertices or small diameter. If $diam(T) = 2$, then $T$ has a
dominating vertex, and $\gamma_{r2}(T)\leq 2$. This beats $n\geq 3$.
If $diam(T) = 3$, then $T$ has a dominating set of size two, which
yields $\gamma_{r2}(T)\leq 4$. This handles the desired bound for
such trees with at least six vertices. When $n=4$ or $n=5$, then $T$
is a spider and the theorem holds by Proposition \ref{spider}.
Moreover, if $T$ is a path of order four, then it achieves this
bound.

Hence we may assume that $diam(T) \geq 4$. Given a subtree $T'$ with
$n'$ vertices, where $n'\geq 3$, the induction hypothesis yields a
2RDF $f'$ of $T'$ with weight at most $3n'/4$. We find such $T'$ and
add a bit more weight to obtain a 2RDF $f$ of $T$. Let $P$ be a
longest path in $T$ chosen to maximize the degree of the penultimate
vertex $v$ on it, and let $u$ be the non-leaf neighbor of $v$.

\medskip
{\it{Case 1.~  $d_T(v) > 2$.}}
\medskip

We obtain $T'$ by deleting $v$ and its leaf neighbors. Define $f$ on $V(T)$ by letting
$f(x)=f'(x)$ except for $f(v)=\{1,2\}$ and $f(x)=\emptyset$ for each
leaf $x$ adjacent to $v$. Since color set $\{1,2\}$ on $v$ takes
care of its neighbors, $f$ is a 2RDF for $T$. Since
$diam(T)\geq 4$, we have $n'\geq 3$, and $w(f)=w(f')+2\leq 3n'/4+2=3(n-3)/4+2<3n/4$.

\medskip
{\it{Case 2.~  $d_T(v)=d_T(u)=2$.}}
\medskip

We obtain $T'$ by deleting $u$ and $v$ and the leaf neighbor $l$ of
$v$. If $n'=2$, then $T$ is a path of order five and has a 2RDF of
weight $3<3n/4$. Otherwise, the induction hypothesis applies. Define
$f$ on $V(T)$ by letting $f(x)=f'(x)$ except for $f(v)=\{1,2\}$ and
$f(u)=f(l)=\emptyset$. Again $f$ is a 2RDF, and the computation
$w(f)<3n/4$ is the same as in Case $1$.

\medskip
{\it{Case 3.~  $d_T(v)=2$ and $d_T(u)> 2$.}}
\medskip

By the choice of path $P$, every penultimate neighbor of $u$ has
degree 2.

\medskip
{\it{Subcase 3.1.~  Every neighbor of $u$ is penultimate or a
leaf.}}
\medskip

Then $diam(T)=4$ and $T$ is a spider. By Proposition \ref{spider},
$\gamma_{r2}(T)<3n/4$, since $T$ is not a path of order four.

\medskip
{\it{Subcase 3.2.~  There exists a neighbor $t$ of $u$ which is
neither penultimate nor a leaf.}}
\medskip

Then $T-tu$ contains two components $T'$ and $T''$ such that $T''$
is a spider containing $u$. Now $|V(T')|=n'\geq 3$ and the induction
hypothesis applies that $\gamma_{r2}(T')\leq 3|V(T')|/4=3n'/4$. By Proposition \ref{spider}, $\gamma_{r2}(T'')\leq 3|V(T'')|/4$. Hence $\gamma_{r2}(T)\leq
\gamma_{r2}(T')+\gamma_{r2}(T'')\leq 3n/4$. \qed

\medskip

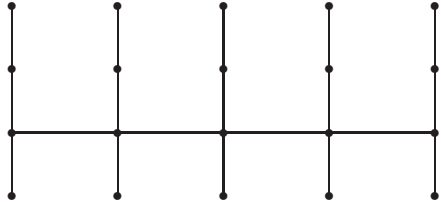
\begin{figure}[h,t]
\setlength{\unitlength}{0.8pt}
\begin{center}
\begin{picture}(40,100)(0,40)

\put(-82,130){\line(0,-1){90}}\put(-32,130){\line(0,-1){90}}\put(18,130){\line(0,-1){90}}\put(68,130){\line(0,-1){90}}\put(118,130){\line(0,-1){90}}

\put(-82,130){\circle*{4}}\put(-32,130){\circle*{4}}\put(18,130){\circle*{4}}\put(68,130){\circle*{4}}\put(118,130){\circle*{4}}

\put(-82,100){\circle*{4}}\put(-32,100){\circle*{4}}\put(18,100){\circle*{4}}\put(68,100){\circle*{4}}\put(118,100){\circle*{4}}

\put(-82,70){\circle*{4}}\put(-32,70){\circle*{4}}\put(18,70){\circle*{4}}\put(68,70){\circle*{4}}\put(118,70){\circle*{4}}

\put(-82,40){\circle*{4}}\put(-32,40){\circle*{4}}\put(18,40){\circle*{4}}\put(68,40){\circle*{4}}\put(118,40){\circle*{4}}

\put(-82,70){\line(1,0){200}}

\end{picture}\vspace{2mm}\caption{The tree $L_5$.}\label{treebound}
\end{center}
\end{figure}

\vspace{2mm}

\medskip

Let $L_k$ consist of the disjoint union of $k$ copies of $P_4$ plus
a path through the center vertices of these copies, as illustrated
in Figure \ref{treebound}. Let $G$  be a graph having an
induced subgraph $P_4$ such that only the center of $P_4$ can be
adjacent to the vertices in $G-P_4$, then every 2RDF of $G$ must have weight at
least 3 on $P_4$. In $L_k$, there are $k$ disjoint $P_4$ of this
form, so $\gamma_{r2}(L_k)\geq 3k=3n/4$. Indeed, we can assemble
such copies of $P_4$ in many ways, and this allows us to
characterize the trees achieving equality in Theorem \ref{tree}.

\begin{thm}\label{treeexact}
Let $T$ be a tree of order $n\geq 3$. Then
$\gamma_{r2}(T)=3n/4$ if and only if $V(T)$ can be partitioned into
sets inducing $P_4$ such that the subgraph induced by the center
vertices of these $P_4$ is connected.
\end{thm}

\proof We have observed that if an induced subgraph $H$ of $G$ is
isomorphic to $P_4$, and its noncenter vertices have no neighbors
outside $H$ in $G$, then every 2RDF of $G$ must have weight at least
3 on $V(H)$. Thus in any tree with the structure described, weight
at least 3 is needed on every $P_4$ in the specified partition. To
show that equality requires this structure, we examine the cases
more closely in the proof of Theorem \ref{tree}. The proof is by
induction on $n$. In the base cases and Cases $1$ and $2$, we
produce a 2RDF with weight less than $3n/4$ except for $P_4$. Define
$u$, $T'$, $T''$, $n'$, $t$ as in the inductive part of Case $3$.
The equality holds only if $n'= n-4$ and $T''$ is a $P_4$ path.
Equality also requires $\gamma_{r2}(T')= 3n'/4$, so by the induction
hypothesis $T'$ has the specified form.

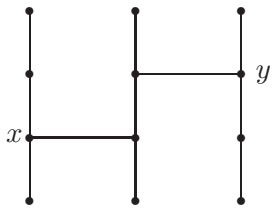
\begin{figure}[h,t]
\setlength{\unitlength}{0.8pt}
\begin{center}
\begin{picture}(40,100)(0,40)

\put(-32,130){\line(0,-1){90}}\put(18,130){\line(0,-1){90}}\put(68,130){\line(0,-1){90}}

\put(-32,130){\circle*{4}}\put(18,130){\circle*{4}}\put(68,130){\circle*{4}}

\put(-32,100){\circle*{4}}\put(18,100){\circle*{4}}\put(68,100){\circle*{4}}

\put(-32,70){\circle*{4}}\put(18,70){\circle*{4}}\put(68,70){\circle*{4}}

\put(-32,40){\circle*{4}}\put(18,40){\circle*{4}}\put(68,40){\circle*{4}}

\put(-32,70){\line(1,0){50}}\put(18,100){\line(1,0){50}}

\put(-43,68){$x$}\put(75,98){$y$}

\end{picture}\vspace{2mm}\caption{A spanning subgraph $H'$ of $T$.}\label{counter}
\end{center}
\end{figure}

Next we show no copy of $P_4$ in $T$ such that both the two
penultimate vertices on $P_4$ with degree at least three in $T$.
Suppose there is a spanning subgraph $H'$ isomorphic to the graph
shown in Figure \ref{counter}, then we give a 2RDF $f$ for $H'$ as
follows:

$$f(v)=\left\{\begin{array}{lllcrrr}
    \{1,2\}  & \mbox{\ \ $v=x$\ or\ $y$},\\[5pt]
    \{1\}  &  \ \ \mbox{$v\notin N[x]\cup N[y]$},\\[5pt]
    \emptyset  & \mbox{\ \ otherwise}.\\[5pt]

\end{array} \right.$$

By Theorem \ref{tree}, $\gamma_{r2}(T)\leq
\gamma_{r2}(H')+\gamma_{r2}(T-H')\leq 8+3(n-12)/4<3n/4$, a
contradiction. \qed

\bigskip

Recall that the {\it corona} $HoK_1$ of a graph $H$ is obtained by attaching one
pendent edge at each vertex of $H$. Since the rainbow domination
number does not increase when edges are added to a graph, we infer
from Theorem \ref{tree} and \ref{treeexact} the following general upper bound.

\begin{cor}\label{connected}
Let $G$ be a connected graph of order $n\geq 3$. Then
$\gamma_{r2}(G)\leq 3n/4$. Moreover, the equality holds if and only
if $G$ is $P_4$ or $C_4oK_1$ or $V(G)$ can be partitioned into $k$
copies of $P_4$ $(k\geq 3)$ and all the copies of $P_4$ can only be
connected by their centers.
\end{cor}

\proof If $G$ has the specified form, then for each copy of $P_4$ in the partition of $V(G)$, every 2RDF of $G$ puts weight at least 3 on it.

Suppose $\gamma_{r2}(G)=3n/4$ and $G$ is not a tree. Since adding
edges can not increase the 2-rainbow domination number, every
spanning tree of $G$ has the form specified in Theorem
\ref{treeexact}. If $n=4$, then $G$ is $P_4$. If $n=8$, then it is
easy to check that the only extremal graph is $C_4oK_1$. If $n\geq
12$, let $T$ be a spanning tree of $G$ has the form specified in
Theorem \ref{treeexact}. $G$ is not a tree, so there exists an edge
$e\in E(G)-E(T)$ such that $T\cup e$ contains a cycle $C$. Without
loss of generality, assume $e$ is not an edge connecting two centers
in $T$. If $C$ contains no edge joining the centers in $T$, i.e., $C$ is formed by some vertices of
a copy $P_4$, then a 2RDF with weight $3n/4-1$ can be found, since we only need to put weight 2 on the vertices of this copy of $P_4$ to take care of this copy of $P_4$. If $C$ goes
through an edges $e'$ joining the centers of two copies of $P_4$ in
$T$, then $\gamma_{r2}(T\cup e-e')< 3n/4$ since tree $T\cup e-e'$ is
not the form specified in Theorem \ref{treeexact}. Hence
$\gamma_{r2}(G)<3n/4$. The proof is complete.\qed

\bigskip
The following result for the 2-rainbow domination number
of paths is given by Bre\v{s}ar and Kraner \v{S}umenjak.

\begin{prop}\label{path}
{\rm{(\cite{bs})}} $\gamma_{r2}(P_n)=\lfloor\frac{n}{2}\rfloor+1$.
\end{prop}

We conclude this subsection with an upper bound in terms of
diameter.

\begin{thm}\label{diameter2}
For any connected graph $G$ on $n$ vertices, $$\gamma_{r2}(G)\leq
n-\lfloor\frac{diam(G)-1}{2}\rfloor.$$ Furthermore, this bound is
sharp.
\end{thm}

\proof Let $P=v_1v_2\cdots v_{diam(G)+1}$ be a diametral path in $G$
and $f$ be a $\gamma_{r2}$-function of $P$. By Proposition
\ref{path}, the weight of $f$ is
$\lfloor\frac{diam(G)+1}{2}\rfloor+1$. Define $g: V(G)\rightarrow
\mathscr{P}(\{1,2\})$ by $g(x)$=$f(x)$ for $x\in V(P)$ and
$g(x)=\{1\}$ for $x\in V(G)-V(P)$. Obviously $g$ is a 2RDF for $G$.
Hence,
$$\gamma_{r2}(G)\leq w(f)+(n-diam(G)-1)=n-\lfloor\frac{diam(G)-1}{2}\rfloor$$
The family of all paths of even order achieves the bound, and the
proof is complete. \qed

\subsection{Lower bounds}

We present a lower bound on the 2-rainbow domination number of a
tree expressed in terms of its domination number, maximum degree,
and the number of its leaves and penultimate vertices. Given a tree, $T$ we denote by $\ell(T)$ the number of leaves in $T$, and by
$p(T)$ the number of penultimate vertices in $T$.

\begin{thm}\label{lower}
For any tree $T$ on at least three vertices, $\gamma_{r2}(T)\geq
\gamma(T)+ \lceil\frac{\ell(T)-p(T)}{\Delta(T)}\rceil$, where $\Delta(T)$ denotes the maximum degree in $T$.
\end{thm}
\proof The proof is by induction on the order of $T$. First we
handle trees with small diameter. If $diam(T)\leq 2$ then
$\gamma(T)=1$, $\gamma_{r2}(T)=2$, and one can easily find that the
required inequality holds. Moreover, we have
$\gamma_{r2}(G)=\gamma(T)+\lceil\frac{\ell(T)-p(T)}{\Delta(T)}\rceil$
precisely when $T$ is isomorphic to $K_{1,r}$ for $r>1$. If
$diam(T)=3$ then another simple analysis shows that the inequality
holds, and the equality is achieved for $DK_{r,s}$ with
$r\geq s\geq 4$ and $DK_{r,1}$ with $r\geq 2$.

Let $T$ be a tree. By the above we may assume that $diam(T)\geq 4$.
Let $P$ be a diametral path with
the leaf $w$ as one of its ends. Suppose $v$ is the neighbor of $w$
and $u$ is the neighbor of $v$ that is not a leaf (hence $u$ also
lies on $P$). Let $L$ denote the vertex set containing $v$ and all leaves
adjacent to $v$ and $F(u)$ be all the possible color sets among all $\gamma_{r2}$-function of $T-L$. Then $\gamma(T-L)\leq \gamma(T)\leq \gamma(T-L)+1$, $\Delta(T)-1\leq \Delta(T-L)\leq \Delta(T)$ and $p(T-L)\leq p(T)\leq p(T-L)+1$.

\medskip
{\it{Case 1.~  $d_T(v)=2$ and $F(u)=\{\{1\},\{2\},\{1, 2\}\}$.}}
\medskip

In this case $\gamma_{r2}(T)=\gamma_{r2}(T-L)+1$. By induction
hypothesis $\gamma_{r2}(T-L)\geq \gamma(T-L)+
\lceil\frac{\ell(T-L)-p(T-L)}{\Delta(T-L)}\rceil$. We finally get
\begin{eqnarray*}
 \gamma_{r2}(T)&=&\gamma_{r2}(T-L)+1 \\
 &\geq & \gamma(T-L)+  \lceil\frac{\ell(T-L)-p(T-L)}{\Delta(T-L)}\rceil
 +1 \\
 &\geq& \gamma(T)+\lceil\frac{\ell(T-L)-p(T-L)}{\Delta(T-L)}\rceil \\
 &\geq & \gamma(T)+  \lceil\frac{\ell(T)-p(T)}{\Delta(T)}\rceil.
\end{eqnarray*}
Since if $p(T-L)=p(T)$, then $\ell(T-L)=\ell(T)$. Otherwise $p(T-L)=p(T)-1$ and $\ell(T-L)=\ell(T)-1$. The last inequality is
obtained.

\medskip
{\it{Case 2.~  $d_T(v)\geq 3$ or $d_T(v)=2$ and $F(u)=\{\emptyset\}$.}}
\medskip

In this case $\gamma_{r2}(T)=\gamma_{r2}(T-L)+2$. Then we get
\begin{eqnarray*}
 \gamma_{r2}(T)&=&\gamma_{r2}(T-L)+2 \\
 &\geq & \gamma(T-L)+  \lceil\frac{\ell(T-L)-p(T-L)}{\Delta(T-L)}\rceil
 +2 \\
 &\geq& \gamma(T)+\lceil\frac{\ell(T-L)-p(T-L)}{\Delta(T-L)}\rceil+1 \\
 &\geq & \gamma(T)+  \lceil\frac{\ell(T)-p(T)}{\Delta(T)}\rceil.
\end{eqnarray*}
In the last inequality we use that the excess of leaves in $T$ with
respect to $T-L$ does not go beyond $\Delta(T)$. \qed

\medskip

In the above proof we mentioned several examples of trees with
diameter at most 3 that achieve the bound in Theorem \ref{lower}. We
pose a characterization of all these extremal graphs as an open problem.

\medskip

Next we give a lower bound of the 2-rainbow domination number of an
arbitrary graph in terms of its diameter.

\begin{thm}\label{diameterupper}
For any connected graph $G$, $\gamma_{r2}(G)\geq
\lceil\frac{2diam(G)+2}{5}\rceil$.
\end{thm}

\proof Let $f=(V_0, V_1^1, V_1^2, V_2)$ be a 2RDF of $G$. Consider
an arbitrary path of length $diam(G)$. This diametral path includes
at most two edges from the induced subgraph $\langle N[v] \rangle_G$
for each vertex $v\in V_1^1\cup V_1^2\cup V_2$. Furthermore, if
vertex $v\in V_0$, then it is adjacent to a vertex with color set
$\{1,2\}$, or adjacent to two different vertices with color set
$\{1\}$ and $\{2\}$, respectively. Hence excluding the edges
mentioned above, the diametral path includes at most
$\min\{|V_1^1|,|V_1^2|\}+|V_2|-1$ other edges joining the
neighborhoods of the vertices of  $V_1^1\cup V_1^2\cup V_2$.
Therefore
\begin{eqnarray*}
 diam(G) & \leq & 2(|V_1^1|+|V_1^2|+|V_2|)+\min\{|V_1^1|,|V_1^2|\}+|V_2|-1 \\
& \leq & 2(|V_1^1|+|V_1^2|+|V_2|)+(|V_1^1|+|V_1^2|)/2+|V_2|-1 \hspace{5em} \\
& = & 5/2(|V_1^1|+|V_1^2|+2|V_2|)-2|V_2|-1 \hspace{5em}\\
& \leq & 5/2\gamma_{r2}(G)-1. \hspace{2em}
\end{eqnarray*}
Then the desired result follows. \qed

\medskip Clearly, the bound of Theorem \ref{diameterupper} is sharp, e.g.
for $G$ isomorphic to $P_3$ or $C_4$.

\vskip 20pt

\noindent \title{\Large\bf Acknowledgments} \maketitle
The authors are indebted to Professor Bostjan Bresar for
his valuable comments and the ideas for this paper.

\end{document}